\documentclass{article}[11pt]
\usepackage{amsthm,amsmath,a4wide}
\usepackage{amssymb}
\input xypic
\usepackage[all,tips]{xy}
\usepackage{hyperref}

\newtheorem{thm}{Theorem}[section]
\newtheorem{lem}[thm]{Lemma}
\newtheorem{prop}[thm]{Proposition}
\newtheorem{cor}[thm]{Corollary}
\newtheorem{conj}[thm]{Conjecture}

\theoremstyle{definition}
\newtheorem{defin}[thm]{Definition}

\theoremstyle{remark}
\newtheorem{remark}[thm]{Remark}
\newtheorem{example}[thm]{Example}

\newcommand{\bth}{\begin{thm}}
\renewcommand{\eth}{\end{thm}}
\newcommand{\bpr}{\begin{prop}}
\newcommand{\epr}{\end{prop}}
\newcommand{\ble}{\begin{lem}}
\newcommand{\ele}{\end{lem}}
\newcommand{\bco}{\begin{cor}}
\newcommand{\eco}{\end{cor}}
\newcommand{\bde}{\begin{defin}}
\newcommand{\ede}{\end{defin}}
\newcommand{\bex}{\begin{example}}
\newcommand{\eex}{\end{example}}
\newcommand{\bre}{\begin{remark}}
\newcommand{\ere}{\end{remark}}
\newcommand{\bcj}{\begin{conj}}
\newcommand{\ecj}{\end{conj}}

\newcommand{\beq}{\begin{equation}}
\newcommand{\eeq}{\end{equation}}

\newcommand{\ot}{{\otimes}}
\newcommand{\op}{{\oplus}}
\newcommand{\lb}{\label}
\newcommand{\nl}{\newline}
\newcommand{\bpf}{\begin{proof}}
\newcommand{\epf}{\end{proof}}

\newcommand{\C}{{\cal C}}
\newcommand{\Z}{{\cal Z}}
\newcommand{\D}{{\cal D}}
\newcommand{\bZ}{{\mathbb Z}}

\newcommand{\e}{{\mathfrak e}}

\newcommand{\Rep}{{\cal R}{\it ep}}
\newcommand{\Vect}{{\cal V}{\it ect}}
\begin{document}
\author{Alexei Davydov}
\title{Unphysical diagonal modular invariants}
\maketitle
\date{}

\begin{center}
Department of Mathematics, Ohio University, Athens, OH 45701, USA
\end{center}

\begin{abstract}
A modular invariant for a chiral conformal field theory is physical if there is a full conformal field theory with the given chiral halves realising the modular invariant.
The easiest modular invariants are the charge conjugation and the diagonal modular invariants. While the charge conjugation modular invariant is always physical there are examples of chiral CFTs for which the diagonal modular invariant is not physical. 
Here we give (in group theoretical terms) a necessary and sufficient 
condition for diagonal modular invariants of $G$-orbifolds of holomorphic conformal field theories to be physical.

Mathematically a physical modular invariant is an invariant of a Lagrangian algebra in the product of (chiral) modular categories.
The chiral modular category of a $G$-orbifold of a holomorphic conformal field theory is the so-called (twisted) Drinfeld centre $\Z(G,\alpha)$ of the finite group $G$.
We show that the diagonal modular invariant for $\Z(G)$ is physical if and only if the group $G$ has a {\em double class inverting} automorphism, that is an automorphism $\phi:G\to G$ with the property that for any commuting $x,y\in G$ there is $g\in G$ such that
$\phi(x) = gx^{-1}g^{-1},\ \phi(y) = gy^{-1}g^{-1}.$

Groups without double class inverting automorphisms are abundant and provide examples of chiral conformal field theories for which the diagonal modular invariant is unphysical. 
\end{abstract}
\tableofcontents
\section{Introduction}

The aim of this note is to construct examples of chiral rational conformal field theories with the same chiral algebras for which the diagonal modular invariant is not physical.
\nl
Recall the state space of a (2-dimesional) conformal field theory comes equipped with amplitudes associated to 
a finite collection of fields inserted in marked points of a Riemann surface \cite{dms,MS}. 
Fields, whose amplitudes depend (anti-)holomorphically on insertion 
points, form what is known as (anti-)chiral algebra of the CFT. The state 
space is naturally a representation of the product of the chiral and anti-chiral 
algebras.
A conformal field theory is {\em rational} if the state space  is a finite sum of tensor 
products of irreducible representations of chiral and anti-chiral algebras.
The matrix of multiplicities of irreducible representations in the decomposition of the state space is called the {\em modular invariant} of the RCFT.
The simplest case (the so-called {\em Cardy case}) is the case of the {\em charge conjugation} modular invariant, which assumes that chiral and anti-chiral algebras coincide. Another very natural case is the {\em diagonal} modular invariant.
\nl 
The name modular invariant comes from the fact that the matrix of multiplicities is invariant with respect to the modular group action on characters. This fact was used in \cite{ciz} to classify modular invariants for affine $sl(2)$ rational conformal field theories.
This paper started an activity aimed at classifying possible modular invariants for various conformal field theories.
It took some time to realise that not all modular invariants correspond to conformal field theories, in other words there are {\em unphysical} modular invariants \cite{fss,ss}.
In the paper \cite{da2} examples of different rational conformal field theories with the same (charge conjugation) modular invariant were constructed.
Thus although being a convenient numerical invariants of a rational conformal field theory modular invariants are far from being complete.
\nl
An adequate description of rational conformal field theories was obtained relatively recently (see \cite{frs} and references therein).
Mathematical axiomatisation of chiral algebras in CFTs is the notion of {\em vertex operator algebra} \cite{Bo,FHL,Ka}. A vertex operator algebra is called {\em rational} if it is a chiral algebra of a RCFT, in particular its category of modules is semi-simple. In this case the category of modules has more structure (see \cite{Hu} and references therein), it is the so-called {\em modular category}. This type of tensor categories was first studied in physics \cite{MS} and then axiomatised mathematically in \cite{Tu}. 
The state space of a RCFT corresponds to a special commutative algebra in the product or categories of modules of chiral and anti-chiral algebras (for details see \cite{frs,HK} and references therein).
This special class of commutative algebras in braided tensor categories is known as {\em Lagrangian} algebras \cite{dmno}. 
The modular invariant is just the class of this Lagrangian algebra in the Grothendieck ring.

Here we characterise Lagrangian algebras with the diagonal modular invariant.
It follows from results of \cite{dno} that such algebras should correspond to braided tensor autoequivalences of the category of representations $\C$ of one of the chiral algebra.
A braided tensor autoequivalence $F:\C\to\C$ corresponds to the diagonal modular invariant if $F(X)\simeq X^*$ for any $X\in\C$. Here $X^*$ is the dual object to $X$. We call such braided tensor autoequivalence {\em dualising}.
We provide examples of modular categories without dualising braided tensor autoequivalences.

Our examples come from permutation orbifolds of {\em holomorphic} conformal field theories (CFTs whose state space is an irreducible module over the chiral algebras).
It is argued in \cite{Ki} (see also \cite{Mu}) that the modular category of the $G$-orbifold of a holomorphic conformal field theory is the so called {\em Drinfeld} (or {\em monoidal}) {\em centre} $\Z(G,\alpha)$, where $\alpha$ is a 3-cocycle of the group $G$.
It is also known that the cocycle $\alpha$ is trivial for {\em permutation orbifolds} (orbifolds where the group $G$ is a subgroup of the symmetric group permuting copies in a tensor power of a holomorphic theory).
The assumption crucial for the arguments of \cite{Ki}  is the existence of twisted sectors.
This assumption is known to be true for permutation orbifolds \cite{bdm}. 
Thus examples we construct correspond to chiral holomorphic orbifolds which do not admit full CFTs with the diagonal modular invariant.
Note that examples (of different nature) of chiral CFTs with unphysical diagonal modular invariant are known \cite{ss}. 

We prove that the category of representations $\Rep(G)$ has a braided dualising autoequivalence if and only if the group $G$ has an automorphism $\phi:G\to G$ with the property
$$\phi(g)\in (g^{-1})^G,\qquad\forall\ g\in G.$$
Here $(g^{-1})^G$ is the conjugacy class of $g^{-1}$. 
We call such automorphisms {\em class-inverting}. 
Using the description of braided monoidal autoequivalences and their action on characters from \cite{da1} we prove that the Drinfeld centre $\Z(G)$ has a braided dualising autoequivalence if and only if the group $G$ has an automorphism $\phi:G\to G$ with the property that for any commuting pair $f,g$ of elements in $G$ there is $h\in G$ such that
$$\phi(f) = hf^{-1}h^{-1},\qquad\phi(g) = hg^{-1}h^{-1}.$$
We call such automorphisms {\em double class-inverting}. 
We proceed by constructing examples of groups without class-inverting automorphisms.


Throughout $k$ denote an algebraically closed field of characteristic zero.
We will work with fusion categories.
A category $\C$ is {\em tensor} over $k$ if it is monoidal and enriched in the category $\Vect_k$ of finite dimensional vector spaces over $k$. That is hom-sets $\C(X,Y)$ for $X,Y\in\C$ are finite dimensional vector spaces over $k$ and the composition and tensor product of morphisms are bilinear maps.
\nl
A tensor category is {\em fusion} if it is semi-simple with finitely many simple objects.
Note that for a semi-simple $\C$ the natural embedding of the Grothendieck group into the dual of the endomorphism algebra of the identity functor
$$K_0(\C)\ \to\ End(Id_\C)^*,\qquad X\mapsto (a\mapsto a_X),\qquad X\in \C,\ a\in End(Id_\C)$$
induces an isomorphism
\beq\lb{gsc}
K_0(\C)\ot_\bZ k\ \to\ End(Id_\C)
\eeq

\section*{Acknowledgment}

The author would like to thank Tom Wolf for suggesting examples of finite simple groups without class-inverting automorphisms (examples \ref{woa} and \ref{ag}). 
Special thanks are to Victor Ostrik for sharing his (and William Kantor's) result on class-inverting automorphisms of groups of odd order (section \ref{goo}), which he used for constructing examples of unphysical modular invariants.

The paper was written while the author was visiting Max Planck Institute for Mathematics (Bonn, Germany) and Macquarie University (Sydney, Australia). The author would like to thank these institutions for hospitality and inspiring atmosphere. 

\section{Group-theoretical braided fusion categories and their dualising braided autoequivalences}

\subsection{Dualising braided autoequivalences}\lb{dba}

Let $\C$ be a pivotal fusion category.
We call a tensor autoequivalence $F:\C\to\C$ {\em dualising}  if $F(X)\simeq X^*$ for any $X\in\C$. Here $X^*$ is the dual object to $X$. 

Note that the composition $F = F'\circ F''$ of any two dualising autoequivalences $F', F''$ has the property $F(X)\simeq X$ for any $X\in\C$.

\ble\lb{dse}
Let $F:\C\to\C$ be an autoequivalence of a semi-simple category $\C$ such that $F(X)\simeq X$ for any $X\in\C$.
Then $F$ is isomorphic to the identity functor.
\ele
\bpf
Choose an isomorphisms $a_Z:F(Z)\to Z$ for a representative $Z\in Irr(\C)$ of every isomorphism class of simple objects.
For an arbitrary $X\in\C$ define $a_X:F(X)\to X$ to fit into a commutative diagram 
$$\xymatrix{F(X) \ar[rr]^{a_X} && X\\ F(\op_{Z\in Irr(\C)} \C(Z,X)Z) \ar[u] \ar[r] & \op_{Z\in Irr(\C)} \C(Z,X)F(Z) \ar[r]^{\op_Z1a_Z} & \op_{Z\in Irr(\C)} \C(Z,X)Z\ar[u]}$$
For $f\in\C(X,Y)$ the commutative diagram
$$\xymatrix{F(X) \ar[rrrr]^{a_X} \ar[ddddd]_{F(f)} &&&& X \ar[ddddd]^f \\
& F(\op_{Z} \C(Z,X)Z) \ar[rd] \ar[lu] && \op_{Z} \C(Z,X)Z \ar[ru] \ar[ddd]_{\op_Z\C(Z,f)1}\\
&& \op_{Z} \C(Z,X)F(Z) \ar[ru]^{\op_Z 1a_Z} \ar[d]_{\op_Z\C(Z,f)}\\
&& \op_{Z} \C(Z,Y)F(Z) \ar[rd]_{\op_Z 1a_Z} \\
& F(\op_{Z} \C(Z,Y)Z) \ar[ld] \ar[ru] && \op_{Z} \C(Z,Y)Z \ar[rd]\\
F(Y) \ar[rrrr]^{a_Y} &&&& Y
}$$
shows that $a:F\to Id_\C$ is natural isomorphism.
\epf

Recall from \cite{da2} that a tensor autoequivalence $F:\C\to\C$ is {\em soft}  if $F$ is isomorphic (as just a functor) to the identity functor.
The following becomes straightforward.
\bco
A dualising tensor autoequivalence $F:\C\to\C$ is unique up to the composition with a soft tensor autoequivalence.
\eco
In other words if a braided fusion category $\C$ has a dualising tensor autoequivalence then the group $Aut^{soft}_\ot(\C)$ of (isomorphism classes of) tensor autoequivalences has an index two extension $\widetilde{Aut^{soft}_\ot}(\C)$ consisting of soft and dualising tensor autoequivalences:
\beq\lb{ite}\xymatrix{1 \ar[r] & Aut^{soft}_\ot(\C) \ar[r] & \widetilde{Aut^{soft}_\ot}(\C) \ar[r] & \bZ/2\bZ \ar[r] & 0}\eeq

The following is straightforward.
\ble\lb{res}
Let $\D$ be a tensor subcategory of a tensor category $\C$.
Then any dualising tensor autoequivalence $F$ of $\C$ preserves $\D$.
\nl
The restriction $F|_\D$ of $F$ to the subcategory $\D$ is a dualising tensor autoequivalence. 
\ele

\bre\lb{bcf}
A tensor autoequivalence $F$ of a fusion category $\C$ is dualising iff its effect on the Grothendieck group $K_0(\C)$ coincides with the dualising map $(\ )^*$. 
\ere

\subsection{Dualising braided autoequivalences of categories of representations}\lb{dbr}

For a finite group $G$ denote by $\Rep(G)$ the category of finite dimensional representations over $k$. 
The category $\Rep(G)$ is a braided (symmetric) tensor category. 
Representation categories are contravariant in $G$: for a group homomorphism $\phi:G\to F$ there is a braided tensor functor
$$\phi^*:\Rep(F)\to\Rep(G)$$
called the {\em inverse image} along $\phi$. More precisely for an $F$ representation $V$ the $G$-action on $\phi^*(V)$ has the form
$g(v) = \phi(g)(v)$, where $g\in G$ and $v\in V$. 
Clearly $\phi^*\circ\psi^* = (\psi\phi)^*$ for group homomorphisms $\phi:G\to H$ and $\psi:H\to F$. 
In particular group automorphisms of $G$ give rise to braided tensor autoequivalences of $\Rep(G)$.
Inner automorphisms (automorphisms of the form $\phi(x) = gxg^{-1}$) are tensor isomorphic to the identity functor.
Thus we have a homomorphism of groups
\beq\lb{arc}Out(G)\to Aut_{br}(\Rep(G))\qquad \phi\mapsto (\phi^{-1})^*\ \eeq
Here $Out(G)$ is the group of outer automorphisms of $G$, that is the quotient of the group of automorphisms $Aut(G)$ by its normal subgroup $Inn(G)$ consisting of inner automorphisms. 
It follows form the Deligne's theorem (on the existence of a fibre functor) \cite{de} that the map \eqref{arc} is an isomorphism.
\bre\lb{gaa}
The results of \cite{da} provide a more elementary proof.
It was proved in \cite{da} that tensor autoequivalences of $\Rep(G)$ correspond to $G$-biGalois algebras.
It is not hard to see that braided tensor autoequivalences of $\Rep(G)$ correspond to commutative $G$-biGalois algebras.
Being semi-simple these algebras have to be isomorphic to the function algebra $k(G)$ with the left and right $G$-actions given by
$$(ga)(x) = a(xg),\qquad (ag)(x) = a(\phi(g)x),\qquad g,x\in G,\ a\in k(G)\ ,$$
where $\phi:G\to G$ is an isomorphism.
\ere

We call a group automorphism $\phi:G\to G$  {\em class-inverting} if for any $g\in G$ there is $h\in G$ such that
$$\phi(g) = hg^{-1}h^{-1}\ .$$

\bpr
The category of representations $\Rep(G)$ of a finite group $G$ has a dualising braided tensor autoequivalence
if and only if the group $G$ has a class-inverting automorphism. 
\epr
\bpf
Since any braided tensor autoequivalence of $\Rep(G)$ has the form $\phi^*$ for an automorphism $\phi$ it is enough to show that $\phi^*$ is dualising if and only if the group $\phi$ is class-inverting.
Note that since $\Rep(G)$ is fusion $\phi^*(V)\simeq V^*$ if and only if the characters of $\phi^*(V)$ and $V^*$ coincide.
The effects on characters
$$\chi_{\phi^*(V)}(g) = \chi_V(\phi(g)),\qquad \chi_{V^*}(g) = \chi_V(g^{-1})$$
imply the proposition.
\epf

It is known (see \cite{da0,da2}) that soft braided tensor autoequivalences of $\Rep(G)$ correspond to so-called {\em class-preserving} automorphisms of $G$, i.e. automorphisms $\phi:G\to G$ such that for every $x\in G$ there is $g\in G$ with $\phi(x) = gxg^{-1}$. 
Denote by $Aut_{cl}(G)$ the group of class-preserving automorphisms of $G$ and by $Out_{2-cl}(G)$ the quotient $Aut_{2-cl}(G)/Inn(G)$ by inner automorphisms. Then the group $Aut_{br}^{soft}(\Rep(G))$ of (isomorphism classes) of soft braided tensor autoequivalences is isomorphic to $Out_{2-cl}(G)$.

The product of two class-inverting automorphisms is class-preserving. 
Thus if $G$ has a class-inverting automorphism the group $Out_{cl}(G)$ has an index two extension
$$\xymatrix{1 \ar[r] & Out_{cl}(G) \ar[r] & \widetilde{Out_{cl}}(G) \ar[r] & \bZ/2\bZ \ar[r] & 0}$$
consisting of outer class-preserving and class-inverting automorphisms.
This is the braided variant of the extension \eqref{ite} for the category $\Rep(G)$.
The index two extension $\widetilde{Aut_{br}^{soft}}(\Rep(G))$ coincides with $\widetilde{Out_{cl}}(G)$.

\subsection{Dualising braided autoequivalences of Drinfeld centres}\lb{dbc}

We call a $G$-action on a $G$-graded vector space $V = \oplus_{g\in G}V_g$ {\em compatible} (with the grading) if $f(V_g) = V_{fgf^{-1}}$.
Let $\Z(G)$ be the category $G$-graded vector spaces with compatible $G$-actions and with morphism being linear maps preserving grading and action. 
\nl
We call  $\Z(G)$ the {\em Drinfeld centre} of a finite group $G$.
Define the tensor product $V\ot U$ of objects $V,U\in\Z(G)$ as the tensor product of $G$-graded vector spaces with the charge conjugation $G$-action.
The category $\Z(G)$ is a tensor category with respect to this tensor product and the trivial associativity constraint.
Moreover $\Z(G)$ is braided with the braiding
$$c_{V,U}(v\otimes u) = f(v)\otimes u,\qquad v\in V_f,\ u\in U\ .$$
The functor 
$$\Rep(G)\to \Z(G)$$
considering a $G$-representation as the trivially $G$-graded (concentrated in the trivial degree) is a braided tensor fully faithful functor (full embedding).

Note that the dual object $V^*=Hom_k(V,k)$ to $V\in\Z(G)$ has the grading $(V^*)_g = (V_{g^{-1}})^*$ and the action 
$$f(l)(v) = l(f^{-1}(v)),\qquad f\in G, l\in V^*, v\in V\ .$$

\bre\lb{chz}
Here we briefly recall the character theory of the category $\Z(G)$. We use slightly different notation comparing e.g. to \cite{dpr}.
The {\em character} of an object $V$ of $\Z(G)$ is a function $\chi_V:\{(f,g)\in G^{\times 2}|\ fg=gf\}\to k$ defined by 
$$\chi_V(f,g) = tr_{V_f}(g),$$ where $tr_{V_f}(g)$ is the trace of the linear operator $g:V_f\to V_f,\quad v\mapsto g(v)$.
It is straightforward that a character is a {\em double class function}, that is 
$$\chi(hfh^{-1},hgh^{-1}) = \chi(f,g),\qquad \forall\ h\in G\ .$$
The character of the tensor product $V\ot U$ has the following expression in terms of the characters of $V,U$:
$$\chi_{U\ot V}(f,g) = \sum_{f_1f_2=f}\chi_U(f_1,g)\chi_V(f_2,g)$$
where the sum is taken over all elements $f_1,f_2$ in the centraliser $C_G(g)$.
\nl
The character of the dual is $\chi_{V^*}(f,g) = \chi_V(f^{-1},g^{-1})$. 
\ere

Recall (e.g. from \cite{da2}) a construction of certain braided tensor autoequivalences $F_{\phi,\gamma}:\Z(G)\to\Z(G)$.
First note that for a group isomorphism $\phi:G\to F$ there is a braided tensor equivalence
$$\phi^*:\Z(F)\to\Z(G)$$
called the {\em inverse image} along $\phi$. For $V\in\Z(F)$ the $G$-action on $\phi^*(V)$ has the form
$g(v) = \phi(g)(v)$, where $g\in G$ and $v\in V$ and the $G$-grading is defined by $V_g = V_{\phi^{-1}(g)}$. 
\nl
For normalised 2-cocycle $\gamma\in Z^2(G,k^*)$ the tensor autoequivalence $F_\gamma:\Z(G)\to\Z(G)$ does not change the $G$-grading of $V\in\Z(G)$ but changes the $G$-action to
$$f*v = \frac{\gamma(f,g)}{\gamma(g,f)}f(v),\qquad v\in V_g.$$
The tensor structure of $F_\gamma$ has the form
$$(F_\gamma)_{U,V}:U\ot V\to U\ot V,\qquad u\ot v\mapsto \gamma(f,g)(u\ot v),\quad u\in U_f,\ v\in V_g\ .$$
The braided autoequivalence $F_{\phi,\gamma}:\Z(G)\to\Z(G)$ is the composition $F_\gamma\circ\phi^*$.

The following statement was proved in \cite[Corollary 6.9]{nr}.
\bpr\lb{prc}
A braided tensor autoequivalence of $\Z(G)$ preserving the subcategory $\Rep(G)\to \Z(G)$ has the form $F_{\phi,\gamma}$ for a group automorphism $\phi:G\to G$ and a 2-cocycle $\gamma\in Z^2(G,k^*)$. 
\epr

We call a group automorphism $\phi:G\to G$  {\em double class-inverting} if for any commuting pair $f,g$ of elements in $G$ there is $h\in G$ such that
$$\phi(f) = hf^{-1}h^{-1},\qquad\phi(g) = hg^{-1}h^{-1}.$$

Here is our main result.
\bth\lb{main}
The Drinfeld centre $\Z(G)$ of a finite group $G$ has a dualising braided autoequivalence
if and only if the group $G$ has a double class-inverting automorphism. 
\eth
\bpf
By lemma \ref{res} a dualising braided tensor autoequivalence $F$ of $\Z(G)$ preserves the subcategory $\Rep(G)\subset \Z(G)$.
By proposition \ref{prc} $F$ is tensor isomorphic to $F_{\phi,\gamma}$ for a group automorphism $\phi:G\to G$ and a 2-cocycle $\gamma\in Z^2(G,k^*)$. 
By the remark \ref{bcf} $F:\Z(G)\to\Z(G)$ is dualising if and only if acts as duality on the Grothendieck ring $K_0(\Z(G))$.
According to the remark \ref{chz} $K_0(\Z(G))$ embeds in the algebra $K_0(\Z(G))\ot_\bZ k$ and the algebra $K_0(\Z(G))\ot_\bZ k$ coincides with the algebra of $k$-valued double class functions (the map from $K_0(\Z(G))\ot_\bZ k$ to the algebra of $k$-valued double class functions given by the character is an embedding and hence an isomorphism, since both $K_0(\Z(G))\ot_\bZ k$ and the algebra of $k$-valued double class functions have the same dimension equal to the number of double conjugacy classes).
\nl
We have the following formula for the character of $F_{\phi,\gamma}(V)$: 
$$\chi_{F_{\phi,\gamma}(V)}(f,g) = \frac{\gamma(f,g)}{\gamma(g,f)}\chi_V(\phi(f),\phi(g))\ .$$
Thus $F_{\phi,\gamma}$ is dualising if and only if the pair $\phi,\gamma$ satisfies the condition 
$$\chi(\phi(f),\phi(g)) = \frac{\gamma(g,f)}{\gamma(f,g)}\chi(f^{-1},g^{-1}),\qquad f,g\in G\ $$
for all double class functions $\chi$.
Taking $\chi$ to be the delta-function on a double conjugacy class we get that $\phi:G\to G$ is a double class-inverting automorphism.

Conversely a double class-inverting automorphism $\phi:G\to G$ a dualising braided autoequivalence induces $\phi^*:\Z(G)\to\Z(G)$.
\epf

\bre
Following \cite{da2} we call a group automorphism $\phi:G\to G$ {\em doubly class-preserving} if it preserves conjugacy classes of commuting pairs of elements of $G$, that is for any $x,y\in G$ such that $xy=yx$ there is $g\in G$ such that $\phi(x) = gxg^{-1},\ \phi(y) = gyg^{-1}$. 
It is straightforward that doubly class-preserving automorphisms are closed under the composition.
Denote by $Aut_{2-cl}(G)$ the group of doubly class-preserving automorphisms of $G$.
Clearly inner automorphisms are doubly class-preserving. The quotient $Aut_{2-cl}(G)/Inn(G)$ is denoted $Out_{2-cl}(G)$.
\nl
Denote by $B(G)$ the subgroup of $H^2(G,k^*)$ consisting of classes of $\gamma$ satisfying 
$$\frac{\gamma(f,g)}{\gamma(g,f)} = 1\qquad\mbox{for any commuting}\ f,g\in G\ .$$
The group $B(G)$ is called {\em Bogomolov multiplier} of $G$.
\nl
Note that the natural action of $Out_{2-cl}(G)$ on $H^2(G,k^*)$ leaves the subgroup $B(G)$ invariant.

It follows from theorem 2.12 of \cite{da2} that the group $Aut^{soft}_{br}(\Z(G))$ of soft braided tensor autoequivalences of $\Z(G)$ is isomorphic to the semi-direct product $Out_{2-cl}(G)\ltimes B(G)$. 
Indeed, by theorem 2.12 a soft braided tensor autoequivalence of $\Z(G)$ has a form $F_{\phi,\gamma}$, where
$$\chi(\phi(f),\phi(g)) = \frac{\gamma(g,f)}{\gamma(f,g)}\chi(f,g),\qquad f,g\in G$$
for all double class functions $\chi$.
Taking $\chi$ to be the delta-function on a double conjugacy class we get that $\phi:G\to G$ is a double class-preserving automorphism.
\ere

It is clear that the product of two double class-inverting automorphisms is double class-preserving. 
Thus if $G$ has a double class-inverting automorphism the group $Out_{2-cl}(G)$ has an index two extension
$$\xymatrix{1 \ar[r] & Out_{2-cl}(G) \ar[r] & \widetilde{Out_{2-cl}}(G) \ar[r] & \bZ/2\bZ \ar[r] & 0}$$
consisting of outer double class-preserving and double class-inverting automorphisms.
Moreover the braided variant of the extension \eqref{ite} fits into a commutative diagram:
$$\xymatrix{B(G) \ar[d] \ar@{=}[r] & B(G)\ar[d] \\
Aut^{soft}_{br}(\Z(G)) \ar[r] \ar[d] & \widetilde{Aut^{soft}_{br}}(\Z(G)) \ar[r]  \ar[d] & \bZ/2\bZ  \ar@{=}[d] \\
Out_{2-cl}(G) \ar[r] & \widetilde{Out_{2-cl}}(G) \ar[r] & \bZ/2\bZ}$$
In other words $\widetilde{Aut^{soft}_{br}}(\Z(G))$ is isomorphic to the semi-direct product $\widetilde{Out_{2-cl}}(G)\ltimes B(G)$.

\subsection{Class-inverting automorphisms of finite groups, examples}\lb{exa}

\subsubsection{Groups without outer automorphisms}\lb{woa}

We start by looking at groups without non trivial outer automorphisms.

A group $G$ is {\em ambivalent} if the identity is a class-inverting automorphisms of $G$, i.e. for any $x\in G$ there is $g\in G$ such that $x^{-1} = gxg^{-1}$.
Thus if $Out(G)=1$ then $G$ has a class-inverting automorphism if and only if $G$ is ambivalent.

The following is straightforward.
\ble
The identity is a class-inverting automorphisms of a finite group $G$ if and only if all the irreducible characters of $G$ are real.
\ele

\bex
Let $G$ be the Mathieu group $M_{11}$ of order $7920=2^4\cdot 3^2\cdot 5\cdot 11$. It has no non trivial outer automorphisms, while there are
elements in $M_{11}$ which are not conjugate to their inverses (there are non-real irreducible characters).
More precisely $M_{11}$ has 144 Sylow 11-subgroups. Thus the normaliser of any of them has order 55, which means that elements of order 11 are not conjugate to their inverses.
\nl
Thus $M_{11}$ does not have class-inverting automorphisms.
\eex

Similarly we call a group $G$ is {\em doubly ambivalent} if the identity is a double class-inverting automorphisms of $G$, i.e. for any commuting pair $x, y\in G$ there is $g\in G$ such that $x^{-1} = gxg^{-1},\ y^{-1} = gyg^{-1}$.
Thus if $Out(G)=1$ then $G$ has a double class-inverting automorphism if and only if $G$ is doubly ambivalent.

The category $\Z(M_{11})$ does not have dualising autoequivalences.

\subsubsection{Abelian groups}

Let $A$ be a finite abelian group. The inverse map $\phi:A\to A,\ \phi(a) = a^{-1}$ is a group automorphism.
\nl
The following is straightforward.
\ble
The inverse map $\phi$ is a (double) class-inverting automorphism of an abelian group $A$. 
\ele

The inverse image $\phi^*$ is a dualising autoequivalence for the categories $\Rep(A),\ \Z(A)$ for an abelian group $A$.

\subsubsection{Groups of odd order}\lb{goo}
The results of this section were communicated to me by Victor Ostrik and are due to 
William Kantor and Victor Ostrik.

\bth Let $G$ be a group of odd order. Then $G$ admits a class reversing
automorphism if and only if $G$ is abelian.
\eth
\bpf 
If $G$ is abelian then $x\mapsto x^{-1}$ is a class reversing automorphism.

To prove the opposite implication we need the following well known 

\ble 
\label{noconj}
Assume that $x\in G$ is conjugated to $x^{-1}$. Then $x=1$.
\ele
\bpf 
Assume $gxg^{-1}=x^{-1}$ for some $g\in G$. 
Then $g^nxg^{-n}=x^{(-1)^n}$ for any $n\in \bZ$. Since the order of $g$
is odd we get $x=x^{-1}$. The result follows.
\epf

Let $t$ be a class reversing automorphism of $G$. Clearly any odd power of $t$ is also class
reversing, so we can assume that the order of $t$ is $2^n, n \in \bZ_{\ge 0}$. If $n=0$ then
$t=\mbox{Id}$ and Lemma \ref{noconj} implies that $G$ is trivial. If $n=1$ then $t$ is an involution
and Lemma \ref{noconj} implies that $t$ has no fixed points. It is well known that this implies
that $G$ is abelian and $t(x)=x^{-1}$ for any $x\in G$, see e.g. Theorem 1.4 \cite[Section 10.1]{Gor}. 

We claim that $n>1$ is impossible. Indeed in this case consider $s=t^{2^{n-1}}$. Then $s$ is an
involution and it maps any conjugacy class to itself. However for $x=a^{-1}s(a)$ we have
$s(x)=s(a)^{-1}a=x^{-1}$. By Lemma \ref{noconj} we get $x=1$ and thus $s(a)=a$ for any
$a\in G$. Thus the order of $t$ is $\le 2^{n-1}$ and we have a contradiction.
\epf

Thus we see that for a non-abelian group $G$ of odd order the category ${\mathcal Z}(G)$ 
has no dualising autoequivalences.

\subsubsection{Symmetric groups}

The following is standard.
\ble\lb{ccsg}
Conjugacy classes of the symmetric group $S_n$ correspond to partitions $n = n_1 + ... + n_r$ into positive integers.
\ele
\bpf
The conjugacy class of $\sigma\in S_n$ is completely determined by the decomposition $X = O_1\cup...\cup O_r$ of the set $X=\{1,...,n\}$  into $\langle \sigma,\pi\rangle$-orbits. Denote by $n_i = |O_i|$ the sizes of the orbits. Then $n = n_1 + ... + n_r$ is a partition.
\epf

\ble
The identity is a class-inverting automorphism of the symmetric group $S_n$. 
\ele
\bpf
By lemma \ref{ccsg} conjugacy classes of $\sigma$ and $\sigma^{-1}$ coincide for all $\sigma\in S_n$.
\epf

Now we examine double class-inverting automorphisms of symmetric groups.

\ble
Double conjugacy classes $(\sigma,\pi)$ of $S_n$ correspond to the following data:
\nl
a partition $n = n_1 + ... + n_r$ into positive integers,
\nl
a decomposition $n_i = n_i(\sigma)n_i(\sigma,\pi)n_i(\pi)$ into positive integers for each $i=1,...,r$,
\nl
an invertible element $a_i\in U(\bZ/n_i(\sigma,\pi)\bZ)$ for each $i=1,...,r$. 
\ele
\bpf
The conjugacy class of a commuting pair $(\sigma,\pi)$ is controlled by the orbit structure of the group $\langle \sigma,\pi\rangle$ on the set $X=\{1,...,n\}$. Let $X = O_1\cup...\cup O_r$ be the decomposition into $\langle \sigma,\pi\rangle$-orbits. Denote by $n_i = |O_i|$ the sizes of the orbits. Clearly $n = n_1 + ... + n_r$.
\nl
Any orbit $O_i$ is identified with the quotient $\langle \sigma,\pi\rangle/K_i$ for a subgroup $K_i$ of finite index. Let $m_i(\sigma)$ and $m_i(\pi)$ be the orders of $\sigma$ and $\pi$ in the quotient. Let $n_i(\sigma,\pi)$ be the order of the intersection $\langle \sigma\rangle\cap\langle\pi\rangle$ in the quotient. Then 
$$|\langle \sigma,\pi\rangle/K_i| = \frac{m_i(\sigma)m_i(\pi)}{n_i(\sigma,\pi)} = n_i(\sigma)n_i(\sigma,\pi)n_i(\pi)\ ,$$
where $n_i(\sigma) = \frac{m_i(\sigma)}{n_i(\sigma,\pi)}$ and $n_i(\pi) = \frac{m_i(\pi)}{n_i(\sigma,\pi)}$
\nl
Finally note that the subgroup $K_i$ has a unique presentation of the form 
$$K_i = \langle\ \sigma^{m_i(\sigma)},\ \pi^{m_i(\pi)},\ \pi^{-n_i(\pi)}\sigma^{n_i(\sigma)a_i}\ \rangle$$
for some $a_i\in U(\bZ/n_i(\sigma,\pi)\bZ)$. 
Indeed the intersection $\langle \sigma\rangle\cap\langle\pi\rangle$ in the quotient can be written as $\langle \sigma^{n_i(\sigma)}\rangle = \langle\pi^{n_i(\pi)}\rangle$. So there is a unique $a_i\in U(\bZ/n_i(\sigma,\pi)\bZ)$ such that $\pi^{n_i(\pi)} = \sigma^{n_i(\sigma)a_i}$ in the quotient.
\epf

\bex
Here we list factorised partitions of $n=3$ and (representatives of) the corresponding double conjugacy classes of $S_3$:
$$\begin{array}{ll}1\cdot 1\cdot 1 + 1\cdot 1\cdot 1 + 1\cdot 1\cdot 1 & (e,e)\\
1\cdot 1\cdot 2 + 1\cdot 1\cdot 1 & (e,(12))\\
1\cdot 2\cdot 1 + 1\cdot 1\cdot 1 & ((12),(12))\\
2\cdot 1\cdot 1 + 1\cdot 1\cdot 1 & ((12),e)\\
1\cdot 1\cdot 3 & (e,(123))\\
1\cdot 3\cdot 1,\quad  a=\pm1 & ((123),(123)^a) \\
3\cdot 1\cdot 1 & ((123),e)
\end{array}$$
\eex

Note that $(\sigma,\pi)$ and the inverse pair $(\sigma^{-1},\pi^{-1})$ has the same factorised partitions
$$n_i(\sigma) = n_i(\sigma^{-1}),\quad n_i(\sigma,\pi) = n_i(\sigma^{-1},\pi^{-1}),\quad n_i(\pi) = n_i(\pi^{-1}),\quad a_i(\sigma,\pi) = a_i(\sigma^{-1},\pi^{-1})\ .$$
Thus we have the following.
\bco
The symmetric group $S_n$ is doubly ambivalent, that is the identity is a double class-inverting automorphism of the symmetric group $S_n$. 
\eco

The identity is a dualising autoequivalence for the categories $\Rep(S_n),\ \Z(S_n)$.

\subsubsection{Alternating groups}\lb{ag}

Denote by $\phi:A_n\to A_n$ the outer automorphism given by $\phi(\sigma) = (12)\sigma(12)^{-1}$. 

We start by recalling the following standard fact.
\ble\lb{ccag}
Let $\sigma\in A_n$ be an even permutation.
The conjugacy class $\sigma^{S_n}$ in $S_n$ is the union of $\sigma^{A_n}\cup\phi(\sigma)^{A_n}$ of two conjugacy classes in $A_n$ iff the centraliser $C_{S_n}(\sigma)$ is contained in $A_n$.
Otherwise $\sigma^{S_n}=\sigma^{A_n}$.
\ele
\bpf
It follows from the coset decomposition $S_n = A_n\cup A_n(12)$ that $\sigma^{S_n} = \sigma^{A_n}\cup\phi(\sigma)^{A_n}$. The classes $\sigma^{A_n}, \phi(\sigma)^{A_n}$ coincide iff $C_{S_n}(\sigma)$ contains an odd permutation.
\epf
\bre\lb{sccc}
Let $\sigma\in A_n$ correspond to a partition $n = n_1 + ... + n_r$. That is $\sigma = \tau_1...\tau_r$, where $\tau_i$ cyclically permute elements of disjoint $n_i$-element sets $X_i$. Denote by $n = m_1 + ... + m_s$ the dual partition, i.e. $m_j = |\{i| n_i=j\}|$. 
The centraliser $C_{S_n}(\sigma)$ coincides with $\prod_j S_{m_j}\ltimes C_j^{\times m_j}$, where copies of cyclic group $C_j$ are generated by $\tau_i$ with $n_i=j$ and $S_{m_j}$ acts on the product $C_j^{\times m_j}$ by permuting the factors.
Thus $C_{S_n}(\sigma)$ is contained in $A_n$ iff $\sigma$ is the products of cycles of distinct odd lengths.
\ere

Here we describe those alternating groups which have a class-inverting automorphism. The argument is borrowed from \cite{wi}. 
First we examine when an element of $A_n$  is conjugate to its inverse.
If the conjugacy class in $S_n$ of $\sigma\in A_n$ does not split in $A_n$, then $\sigma$ is conjugate to $\sigma^{-1}$ in $A_n$ (because it is conjugate in $S_n$). Thus, it suffices to check whether an element whose conjugacy class does split inside $A_n$, is conjugate to its inverse. By remark \ref{sccc}, it suffices to look at those even permutations that arise as products of cycles of distinct odd lengths. 
\ble
A product of cycles of distinct odd lengths $n_1,...,n_r$ is conjugate to its inverse if and only if $\sum_{i=1}^r\frac{n_i-1}{2}$ is even. Equivalently, it is conjugate to its inverse if and only if the number of $n_i$ that are congruent to $3$ modulo $4$ is even.
\ele
\bpf
In order to determine whether such a product of cycles $\sigma$ is conjugate to its inverse in $A_n$, it suffices to find a permutation $\pi\in S_n$ that conjugates $\sigma$ to its inverse. Then $\sigma$ is conjugate to its inverse in $A_n$ iff $\pi$ is even.

Note that for a cycle of odd length $\tau=(1...l)$, the product of $(l-1)/2$ transpositions $(i,l+2-i)$ conjugates $\tau$ to its inverse. Thus, a product of cycles of odd lengths $n_1,...,n_r$ is conjugate to its inverse in $S_n$ by a product of $\sum_{i=1}^r\frac{n_i-1}{2}$ transpositions.
\epf

\bpr
The alternating group $A_n$ has a class-inverting automorphism if and only if $n= 1,2,3,4,5,6,7,8,10,12,14$.
\nl
For $n=1,2,5,6,10,14$ the identity is class-inverting.
\nl 
For $n=3,4,7,8,12$ the automorphism $\phi$ is class-inverting.
\epr
\bpf
For given $n$ one of the following occurs:
\nl
(1). For all conjugacy classes that split in $A_n$, the number of cycle sizes that are congruent to 3 mod 4 is even: In this case, every element is conjugate to its inverse, the group is ambivalent, and the identity map is a class-inverting automorphism.
\nl
(2). For all conjugacy classes that split in $A_n$, the number of cycle sizes that are congruent to 3 mod 4 is odd: In this case, no element whose conjugacy class splits is conjugate to its inverse in $A_n$. Thus $\phi$ sends every such element into the conjugacy class of its inverse. For $\sigma$ whose conjugacy class does not split in $A_n$, $\phi(\sigma)$ is still in the conjugacy class of $\sigma^{-1}$. Thus $\phi$ is a class-inverting automorphism.
\nl
(3). Of the conjugacy classes that split, there is at least one with an even number of cycle sizes that are congruent to 3 mod 4 and at least one with an odd number of cycle sizes that are congruent to 3 mod 4. In this case, neither conjugation by an even permutation nor conjugation by an odd permutation is class-inverting. Thus, no conjugation by an element in the symmetric group is class-inverting. But by fact (2), unless $n=6$, these are all the automorphisms, so the group has no class-inverting automorphism.

If $n=9+4k,\ k\geq 0$ then case (3) holds. Indeed $n = (n-4)+3+1$ is a partition with distinct odd parts, an odd number of which is 3 mod 4. At the same time, $n=n$ is a partition with distinct odd parts, an even number of which are 3 mod 4.
\nl
If $n=11+4k,\ k\geq 0$ then case (3) holds. In this case $n = (n-4)+3+1$ is a partition with distinct odd parts, an even number of which is 3 mod 4. On the other hand, $n=n$ is a partition with distinct odd parts, an odd number of which are 3 mod 4.
\nl
If $n=16+4k,\ k\geq 0$ then case (3) holds. Indeed then $n = (n-9)+5+3+1$ gives a partition with distinct odd parts, an even number of which are 3 mod 4. On the other hand, $n=(n-5)+5$ is a partition with distinct odd parts, an odd number of which are 3 mod 4.
\nl
If $n=18+4k,\ k\geq 0$ then case (3) holds. Then $n = (n-9)+5+3+1$, a partition with distinct odd parts, an odd number of which are 3 mod 4. At the same time, $n=(n-5)+5$ is a partition with distinct odd parts, an even number of which are 3 mod 4.
\nl
The only possibilities left are 1,2,3,4,5,6,7,8,10,12,14. Case (1) applies to 1,2,5,6,10,14 (for these $n$ the alternating group $A_n$ is ambivalent) and case (2) applies to 3,4,7,8,12.
\epf

In particular, for $n\geq 15$ the alternating group $A_n$ does not have a class-inverting automorphism and the categories $\Rep(A_n),\ \Z(A_n)$ do not possess a dualising autoequivalence.

\section{Applications}\lb{app}

Here we sketch a construction of chiral conformal field theories for which the diagonal modular invariant is not realisable. 

\subsection{Lagrangian algebras of braided autoequivalences}

Let $F:\C\to\C$ be a braided tensor autoequivalence of a braided fusion category $\C$.
Consider an object 
$$Z(F) = \bigoplus_{X\in Irr(\C)}X\boxtimes F(X)^*\in\C\boxtimes\overline\C\ ,$$
where the sum is taken over isomorphism classes of simple objects of $\C$. 
It is known (see e.g. \cite{da2} for details) that this object has a structure $\mu:Z(F)\ot Z(F)\to Z(F)$ a commutative algebra in $\C\boxtimes\overline\C$. Moreover it follows from the results of \cite[section 3.2]{dno} that $Z(F)$ is a Lagrangian algebra in $\C\boxtimes\overline\C$. 
\nl
Note that the algebras $Z(F), Z(F')$ are isomorphic if and only if the autoequivalences $F, F'$ are tensor isomorphic.

For a Lagrangian algebra $Z\in\C\boxtimes\overline\C$ the class $[Z]$ in the Grothendieck ring $K_0(\C\boxtimes\overline\C)\simeq K_0(\C)\times_\bZ K_0(\C)$ written in the basis of classes of simple objects of $\C$
$$[Z] = \sum_{\chi,\xi}Z_{\chi\xi}(\chi\boxtimes\overline\xi^*)$$
gives rise to a non-negative integer matrix $(Z_{\chi\xi})$ called the {\em modular invariant} of $Z$.
\nl
We say that $Z$ has the {\em diagonal modular invariant} if its matrix is the identity $Z_{\chi\xi} = \delta_{\chi\xi}$.

The next lemma follows from the results of \cite[section 3.2]{dno}.
\ble
A Lagrangian algebra $Z\in\C\boxtimes\overline\C$ has the diagonal modular invariant if and only if $Z \simeq Z(F)$ for a dualising braided tensor autoequivalence $F:\C\to\C$. 
\ele

\subsection{Holomorphic permutation orbifolds}

Let $V$ be a holomorphic vertex operator algebra (for example $V = \e_{8,1}$). 
\nl
Let $G\subset S_n$ be a subgroup of the permutation group.
The vertex operator subalgebra $(V^{\ot n})^G$ of invariants is called the {\em chiral permutation orbifold} of $V$.
\nl
According to \cite{Ki} its category of representation is $\Rep((V^{\ot n})^G) = \Z(G)$ (subject to the rationality of $(V^{\ot n})^G$). 
\nl
According to \cite{frs,HK} a Lagrangian algebra $Z\in \Z(G)\boxtimes\overline\Z(G)$ gives rise to a rational conformal field theory with the left (right) chiral algebras $(V^{\ot n})^G$ and the modular invariant $[Z]$.
In particular rational conformal field theories with the diagonal modular invariant correspond to Lagrangian algebras $Z=Z(F)$ for a dualising braided tensor autoequivalence $F:\Z(G)\to\Z(G)$.
The examples of groups $G$ without class-inverting automorphisms from section \ref{exa} provide examples of chiral rational conformal field theories for which the diagonal modular invariant is non-realisable.

\end{document}